\documentclass{elsart}    

\usepackage{epsfig}
\usepackage{latexsym}
\usepackage{times}
\usepackage{amsfonts}

\newcommand{\cal}[1]{\mbox{\it{#1}}}

\begin{document}

\begin{frontmatter}

\title{Batch and median neural gas}
\author{Marie Cottrell}
\address{SAMOS-MATISSE, Universit\'e Paris I, Paris, France}
\author{Barbara Hammer}
\footnote{corresponding author:
Barbara Hammer, Institute of Computer Science,
Clausthal University of Technology, Julius-Albert-Str.~4, D-38678 Clausthal-Zellerfeld, Germany, hammer@in.tu-clausthal.de}
\address{Institute of Computer Science, Clausthal University of
Technology, Clausthal-Zellerfeld, Germany}
\ead{hammer@in.tu-clausthal.de}
\author{Alexander Hasenfu\ss}
\address{Institute of Computer Science, Clausthal University of
Technology, Clausthal-Zellerfeld, Germany}
\author{Thomas Villmann}  
\address{Clinic for Psychotherapy, Universit\"at Leipzig, Leipzig, Germany}
\newpage

\begin{abstract}
Neural Gas (NG) constitutes a very robust clustering
algorithm given euclidian data which does not suffer from the
problem of local minima like simple vector quantization,
or topological restrictions like the self-organizing map.
Based on the cost function of NG, we introduce a batch variant
of NG which shows much faster convergence and which can be interpreted as
an optimization of the cost function by the Newton method.
This formulation has the additional benefit that,
based on the notion of the generalized median in analogy to Median SOM, a
variant for non-vectorial proximity data can be introduced.
We prove convergence of batch and median versions of NG, SOM, and k-means in
a unified formulation, and we investigate
the behavior of the algorithms in several experiments.
\end{abstract}

\begin{keyword}
Neural gas \sep batch algorithm \sep proximity data \sep median-clustering \sep convergence
\end{keyword}
\end{frontmatter}

\section{Introduction}  
Clustering constitutes a fundamental problem in various
areas of applications such as pattern recognition, image processing, data
mining, data compression, or machine learning \cite{murty}. 
The goal of clustering is grouping given training data into
classes of similar objects such that 
data points with similar semantical meaning are linked together.
Clustering methods differ in  various aspects including
the assignment of data points to classes which might be crisp or fuzzy,
the arrangement of clusters which might be flat or hierarchical, or
the representation of clusters which might be represented by the
collection of data points assigned to a given class or by few prototypical
vectors.
In this article, we are interested in neural clustering algorithms
which deal with crisp assignments and representation of 
clusters by neurons or prototypes.

Popular neural algorithms representing data by a small number of typical
prototypes include
k-means, the self-organizing map (SOM), neural gas (NG), and
alternatives \cite{duda,ripley}.
Depending on the task and model at hand,
these methods can be used for data compression, data mining and visualization,
nonlinear projection and interpolation,
or preprocessing for supervised learning. 
K-means clustering 
directly aims at a minimization of the quantization error \cite{bottou}.
However, its update scheme is local, therefore it easily
gets stuck in local optima.
Neighborhood cooperation as for SOM and NG offers one biologically plausible
solution. Apart from a reduction of
the influence of initialization,
additional semantical insight is gained: browsing within the
map and, if a prior low dimensional lattice is chosen, 
data visualization become possible.
However, a fixed prior lattice as chosen in SOM might be
suboptimal for a given task depending on the data topology
and topological mismatches can easily occur
\cite{villy}.
SOM does not possess a cost function in the continuous case,
and the mathematical analysis is quite difficult unless
variations
of the original learning rule are considered
for which cost functions can be found
\cite{marie,heskes}. 
NG optimizes a cost function which, as a limit
case, yields the
quantization error \cite{ng}.
Thereby, a data optimum (irregular)
lattice can be determined automatically during training which
perfectly mirrors the data topology and which 
allows to browse within the result \cite{trn}. This yields
very robust clustering behavior. Due to the
potentially irregular lattice, visualization 
requires additional projection methods.

These neural algorithms (or a variation thereof for SOM) optimize
some form of cost function connected to the
quantization error of the data set.
There exist mainly two different optimization schemes for these objectives:
online variants, which adapt the prototypes after each pattern,
and batch variants which adapt the prototypes according to all patterns
at once. Batch approaches are usually much faster in particular
for high dimensional vectors, since
only one adaptation is necessary in each cycle
and convergence can usually 
be observed after few steps.
However, the problem of local optima for k-means
remains in the batch variant.
For SOM, topological ordering might be very difficult to
achieve since, at the beginning, ordering does usually not
exist and, once settled in a topological mismatch, the
topology  can hardly be corrected.
The problem of topological mismatches is much more pronounced in Batch SOM than
in online SOM
as shown in \cite{fort}
such that a good (and possibly costly) initialization is essential for the success.
However, due to their efficiency,
batch variants are often chosen
for SOM or k-means if data are available a priori, whereby the existence
of local optima and topological mismatches might cause severe problems.
For NG, some variants of 
batch adaptation schemes occur at singular points in the literature
\cite{zhong},
however, so far, no NG-batch scheme has been 
explicitely derived from the NG cost function
together with a proof of the convergence of the algorithm.
In this article, 
we put the cost functions of NG, (modified) SOM, and k-means into
a uniform notation and derive batch versions thereof together with
a proof for convergence.
In addition, we relate Batch NG to an optimization
by means of the Newton method, and we compare the methods on different
representative clustering problems.

In a variety of tasks such as classification of protein structures,
text documents, surveys, or biological signals, an explicit 
metric vector space
such as the standard euclidian vector space is not available, rather
discrete transformations of data e.g.\ the edit
distance or pairwise proximities
are available \cite{marienn,graepel,seo}.
In such cases, a clustering method which does not
rely on a vector space has to be applied such as
spectral clustering \cite{belkin}.
Several alternatives to SOM have been proposed which can
deal with more general, mostly discrete data
\cite{marienn,graepel,seo}.
The article \cite{kohonen_largescale} proposes a particularly
simple and intuitive possibility for clustering proximity data: the mean value of the Batch SOM
is substituted by the generalized median resulting in Median SOM,
a prototype-based neural network in which prototypes location are
adapted within the data space by batch computations.
Naturally, the same idea can be transferred to Batch NG and k-means
as we will demonstrate in this contribution.
As for the euclidian versions, it can be shown that the median variants of 
SOM, NG, and k-means converge after a finite number of adaptation steps.
Thus, the formulation of neural clustering schemes by means of
batch adaptation opens the way towards the important 
field of clustering complex data structures for which pairwise
proximities or a kernel matrix constitute the interface to the neural
clustering method.

\section{Neural gas}
Assume data points 
$\vec x\in\mathbb{R}^m$ are distributed according to
an underlying distribution
$P$, the goal of NG as introduced in
\cite{ng} is to find prototype locations $\vec w^i\in\mathbb{R}^m$,
$i=1,\ldots,n$, such that these prototypes represent the distribution $P$ as
accurately as possible, minimizing the cost function
$$E_{\mathrm{NG}}(\vec w)=\frac{1}{2C(\lambda)}\sum_{i=1}^n
\int h_{\lambda}(k_i(\vec x,\vec w))\cdot d(\vec x,\vec w^i)\,P(d\vec x)$$
where
$$d(\vec x,\vec y) = (\vec x-\vec y)^2$$ denotes the squared euclidian distance,
$$k_i(\vec x,\vec w^i)=|\{\vec w^j\:|\:d(\vec x,\vec w^j)<d(\vec x,\vec w^i)\}|$$
is the rank of the prototypes sorted according
to the distances, $h_{\lambda}(t)=\exp(-t/\lambda)$ is a Gaussian shaped curve
with neighborhood range $\lambda>0$, and
$C(\lambda)$ is the constant $\sum_{i=1}^nh_{\lambda}(k_i)$.
The learning rule consists of a stochastic gradient descent,
yielding
$$\Delta \vec w^i = \epsilon\cdot h_{\lambda}(k_i(\vec x^j,\vec w))\cdot(\vec x^j-\vec w^i)$$
for all prototypes $\vec w^i$
given a data point $\vec x^j$. Thereby, the neighborhood range
$\lambda$ is decreased during training 
to ensure independence of initialization
at the beginning of training 
and optimization of the quantization error
in the final stages.
As pointed out in \cite{trn}, the result can be associated with
a data optimum lattice such that browsing within the data space
constitutes an additional feature of the solution.

Due to its simple adaptation rule, the independence of a prior lattice, and
the independence of initialization because of the integrated neighborhood 
cooperation,
NG is a simple and highly effective algorithm for data clustering.
Popular alternative clustering algorithms are offered by the
SOM as introduced by Kohonen \cite{kohonen} and
k-means clustering \cite{duda}.

SOM uses the adaptation strength
$h_{\lambda}(\mathit{nd}(I(\vec x^j),i))$ instead of
$h_{\lambda}(k_i(\vec x^j,\vec w))$, $I(\vec x^j)$ denoting the
index of the
closest prototype, the winner, for $\vec x^j$, and $\mathit{nd}$
a priorly chosen, often two-dimensional neighborhood structure of
the neurons. A low-dimensional lattice offers
the possibility to easily visualize data. However, if the primary goal is
clustering, a fixed topology puts restrictions
on the map and topology preservation often
cannot be achieved
\cite{villy}.
SOM does not possess a cost function in the continuous case
and its mathematical investigation is difficult \cite{marie}.
However, if the winner is chosen as the neuron $i$ with
minimum averaged distance 
$\sum_{l=1}^n h_{\lambda}(\mathit{nd}(i,l)) d(\vec x^j,\vec w^l),$
it optimizes the cost
$$E_{\mathrm{SOM}}(\vec w) \sim \sum_{i=1}^n\int 
\chi_{I^*(\vec x)}(i)\cdot 
\sum_{l=1}^n h_{\lambda}(\mathit{nd}(i,l))\cdot d(\vec x,\vec w^l)\,P(d\vec x)$$
as pointed out by {Heskes} \cite{heskes}.
Here, $I^*(\vec x)$ denotes the winner index according to
the averaged distance and $\chi_j(i)$ is the characteristic function of $j$.

K-means clustering adapts only the winner in each step, thus it optimizes
the standard quantization error
$$E_{\mathrm{kmeans}}(\vec w) \sim \sum_{i=1}^n\int 
\chi_{I(\vec x)}(i)\cdot d(\vec x,\vec w^i)\, P(d\vec x)$$
where $I(\vec x)$ denotes the winner index for $\vec x$
in the classical sense.
Unlike SOM and NG, k-means is very sensitive to initialization of
the prototypes since it adapts the prototypes only locally according
to their nearest data points.
An initialization of the prototypes within the data points is
therefore mandatory.

\subsection{Batch clustering}
If training data $\vec x^1$, \ldots, $\vec x^p$ are given priorly, 
fast alternative  batch training schemes exist for both,
k-means and SOM.
Starting from random positions of the prototypes,
batch learning iteratively performs the following two steps
until convergence
\begin{enumerate}
\item[(1)]
determine the winner $I(\vec x^i)$ resp.\ $I^*(\vec x^i)$
for each data point $\vec x^i$,
\item[(2)]
determine new prototypes as
$$\vec w^i={\sum_{j\:|\: I(\vec x^j)=i} \vec x^j}/
{| \{j\:|\: I(\vec x^j)=i \}|}$$ for k-means
and
$$\vec w^i = {\sum_{j=1}^ph_{\lambda}(\mathit{nd}(I^*(\vec x^j),i))\cdot\vec x^j}/
{\sum_{j=1}^ph_{\lambda}(\mathit{nd}(I^*(\vec x^j),i))}$$
for SOM.
\end{enumerate}
Thereby, the neighborhood cooperation is annealed for SOM 
in the same way as in the online
case.

It has been shown in \cite{bottou,cheng}
that Batch k-means and Batch SOM
optimize the same cost functions 
as their online variants,
whereby the modified winner notation as proposed by Heskes is used for SOM.
In addition, as pointed out in \cite{heskes},
this formulation allows to link the models to statistical
formulations and it can be interpreted as a
limit case of EM optimization schemes for appropriate mixture models.

Often, batch training converges after only few (10-100) cycles such
that this training mode offers considerable speedup in comparison
to the online variants:
adaptation of the (possibly high dimensional) prototypes is only
necessary after the presentation of all training patterns instead 
of each single one.

Here, we introduce Batch NG.
As for SOM and k-means, it
can be derived from the cost function of NG, which,
for discrete data $\vec x^1$, \ldots, $\vec x^p$, reads as
$$E_{\mathrm{NG}}(\vec w) \sim 
\sum_{i=1}^n\sum_{j=1}^p h_{\lambda}(k_i(\vec x^j,\vec w))\cdot d(\vec x^j,\vec w^i)\,,$$
$d$ being the standard euclidian metric.
For the batch algorithm, the quantities 
$k_{ij}:=k_i(\vec x^j,\vec w)$ are treated as hidden variables
with the constraint that the
values $k_{ij}$ ($i=1,\ldots,n$) constitute
a permutation of $\{0,\ldots,n-1\}$ for each point $\vec x^j$.
$E_{\mathrm{NG}}$ is interpreted as a function depending on
$\vec w$ and $k_{ij}$
which is optimized in turn with respect to
the hidden variables $k_{ij}$ and with respect to the
prototypes $\vec w^i$,
yielding the two adaptation steps of Batch NG
which are iterated until convergence:
\begin{enumerate}
\item[{\bf (1)}]
determine 
$$k_{ij}=k_i(\vec x^j,\vec w)=|\{\vec w^l\:|\:d(\vec x^j,\vec w^l)<d(\vec x^j,\vec w^i)\}|$$
as the rank of prototype $\vec w^i$,
\item[{\bf (2)}]
based on the hidden variables $k_{ij}$, set
$$\vec w^i=\frac{\sum_{j=1}^p h_{\lambda}(k_{ij})\cdot\vec x^j}{\sum_{j=1}^p h_{\lambda}(k_{ij})}\,.$$
\end{enumerate}

As for Batch SOM and k-means, adaptation takes place only after the
presentation of all patterns with a step size which is optimized
by means of the partial cost function.
Only few adaptation steps are usually necessary due to the 
fact that Batch NG can be interpreted as Newton optimization
method which takes second order information
into account whereas online NG is given by a simple
stochastic gradient descent.

To show this claim, we
formulate the Batch NG update in the form
$$\Delta \vec w^i=\frac{\sum_{j=1}^p h_{\lambda}(k_{ij})\cdot\left(\vec x^j-\vec w^i\right)}{\sum_{j=1}^p h_{\lambda}(k_{i
j})}\,.$$
Newton's method for an optimization of $E_{\mathrm{NG}}$ yields the formula
$$\triangle \vec w^i=-J(\vec w^i)\cdot H^{-1}(\vec w^i)\,,$$
where $J$ denotes the Jacobian of $E_{\mathrm{NG}}$
and $H$ the Hessian matrix.
Since $k_{ij}$ is locally constant, we get up to sets of
measure zero
$$J(\vec w^i)=2\cdot \sum_{j=1}^ph_{\lambda}(k_{ij})\cdot(\vec w^i-\vec x^j)$$
and the Hessian matrix equals a diagonal matrix with entries
$$2\cdot \sum_{j=1}^ph_{\lambda}(k_{ij})\,.$$
The inverse gives the scaling factor
of the Batch NG adaptation, i.e.\ Batch NG equals
Newton's method for the optimization of $E_{\mathrm{NG}}$.

\subsection{Median clustering}
Before turning to the problem of
clustering proximity data,
we formulate Batch NG, SOM, and k-means
within a common cost function.
In the discrete setting,
these three models optimize a cost function of the form
$$E := \sum_{i=1}^n\sum_{j=1}^p f_1(k_{ij}(\vec w))\cdot f_2^{ij}(\vec w)$$
where
$f_1(k_{ij}(\vec w))$ is the characteristic function of the winner, i.e.\
$\chi_{I(\vec x^j)}(i)$ resp.\ $\chi_{I^*(\vec x^j)}(i)$, for
k-means and SOM, and it is 
$h_{\lambda}(k_i(\vec x^j,\vec w))$ for neural gas.
$f_2^{ij}(\vec w)$ equals the distance
$d(\vec x^i,\vec w^j)$ for k-means and NG,
and it is the averaged distance 
$\sum_{l=1}^nh_{\lambda}(\mathit{nd}(i,l))\cdot d(\vec x^j,\vec w^l)$ for SOM.
The batch algorithms optimize 
$E$ with
respect to $k_{ij}$ in step {\bf (1)} assuming fixed $\vec w$.
Thereby, for each $j$, the vector $k_{ij}$ ($i=1,\ldots,n$) is restricted
to a vector with exactly one entry $1$ and $0$, otherwise,
for k-means and SOM. It is restricted to a permutation
of $\{0,\ldots,n-1\}$ for NG.
Thus, the elements $k_{ij}$ come from a discrete set
which we denote by $K$.
In step {\bf (2)}, $E$ is optimized
with respect to $\vec w^j$ assuming
fixed $k_{ij}$. The update formulas
as introduced above can be derived by taking the derivative of $f_2^{ij}$
with respect to $\vec w$.

For proximity data $\vec x^1$, \ldots, $\vec x^p$, only the distance matrix
$d_{ij}:=d(\vec x^i,\vec x^j)$ is available but
data are not embedded in a vector space and
no continuous adaptation is possible, nor does
the derivative of the distance function $d$ exist.
A solution to tackle this setting
with SOM-like learning algorithms
proposed by Kohonen is offered by the Median SOM: it
is based on the notion of the generalized median \cite{kohonen_largescale}.
Prototypes are chosen from the \emph{discrete} set 
given by the training points
$X=\{\vec x^1,\ldots,\vec x^p\}$
in an optimum way. 
In mathematical terms, $E$ is optimized
within the set $X^n$ given by the training data
instead of $(\mathbb{R}^m)^n$.
This leads to the choice of $\vec w^i$
as 
$$\vec w^i=\vec x^{l} \quad \mbox{where}\quad
\quad l=\mbox{argmin}_{l'}\sum_{j=1}^p 
h_{\lambda}(\mathit{nd}(I^*(\vec x^j),i))\cdot d(\vec x^j,\vec x^{l'})$$
in step {\bf (2)}.
In \cite{kohonen_largescale}, Kohonen considers only the data points 
mapped to a
neighborhood of neuron $i$ as potential candidates for $\vec w^i$ and, in
addition, reduces the above sum to points mapped into a neighborhood of $i$.
For small neighborhood range and approximately ordered maps,
this does not change the result but
considerably speeds up the computation.

The same principle can be applied to k-means and Batch NG.
In step {\bf (2)}, instead of
taking the vectors in $(\mathbb{R}^m)^n$ which minimize $E$,
prototype $i$ is chosen as the data point in $X$ with 
$$\vec w^i=\vec x^l\quad\mbox{where}\quad
l=\mbox{argmin}_{l'}\sum_{j=1}^p\chi_{I(\vec x^j)}(l)\cdot d(\vec x^j,\vec x^{l'})$$
assuming fixed $\chi_{I(\vec x^j)}(l)$ for Median k-means and
$$\vec w^i=\vec x^l\quad\mbox{where}\quad
l=\mbox{argmin}_{l'}\sum_{j=1}^ph_{\lambda}(k_{ij})\cdot 
d(\vec x^j,\vec x^{l'})$$
assuming fixed
$k_{ij}=k_i(\vec x^j,\vec w)$ for Median NG.
For roughly ordered maps, a restriction of potential candidates $\vec x^l$ to
data points mapped to a neighborhood of $i$ can speed up training
as for Median SOM.

Obviously, a direct implementation of the new prototype locations requires
time ${\cal O}(p^2n)$, $p$ being the number of patterns and $n$
being the number of neurons, since
for every prototype and every possible prototype location in $X$
a sum of $p$ terms needs to be evaluated.
Hence, an implementation of Median NG requires the
complexity 
${\cal O}(p^2n+pn\log n)$ for each cycle, including the computation of $k_{ij}$
for every $i$ and $j$.
For Median SOM, 
a possibility to speed up training has recently been presented in
\cite{fabrice} which yields an exact computation with costs
only ${\cal O}(p^2+pn^2)$ instead of ${\cal O}(p^2n)$ for the sum. 
Unfortunately, the same technique does not improve the 
complexity of NG. However, further heuristic possibilities to speed-up 
median-training are discussed in \cite{fabrice} which can be transferred to
Median NG.
In particular, the fact that data and prototype assignments are in large
parts identical for consecutive runs at late stages of training
and a restriction to candidate median points in the 
neighborhood of the previous one
allows a reuse of already computed values and a considerable speedup.

\subsection{Convergence}
All batch algorithms optimize $E=E(\vec w)$ by
consecutive optimization of the hidden variables $k_{ij}(\vec w)$ and $\vec w$.
We can assume
that, for given $\vec w$, the values $k_{ij}$ determined
by the above algorithms are unique, introducing 
some order in case of ties.
Note that the values $k_{ij}$ come from a
discrete set $K$.
If the values $k_{ij}$ are fixed, the choice of the optimum $\vec w$ is unique
in the algorithms
for the continuous case, as is obvious from the formulas given above,
and we can assume uniqueness for the median variants by
introducing an order.
Consider the function
$$Q(\vec w',\vec w)=\sum_{i=1}^n\sum_{j=1}^pf_1(k_{ij}(\vec w))\cdot f_2^{ij}(\vec w')\,.$$
Note that $E(\vec w)=Q(\vec w,\vec w)$.
Assume prototypes $\vec w$ are given, and 
new prototypes $\vec w'$ are computed based on
$k_{ij}(\vec w)$ using one of the above batch or median algorithms. 
It holds
$E(\vec w')=Q(\vec w',\vec w') \le Q(\vec w',\vec w)$
because 
$k_{ij}(\vec w')$ are optimum assignments for
$k_{ij}$ in $E$, given $\vec w'$.
In addition,
$Q(\vec w',\vec w)\le Q(\vec w,\vec w)=E(\vec w)$
because $\vec w'$ are optimum assignments
of the prototypes given $k_{ij}(\vec w)$.
Thus, $E(\vec w')-E(\vec w)=
E(\vec w') - Q(\vec w',\vec w) +Q (\vec w',\vec w) -E(\vec w) \le 0$,
i.e., in each step of the algorithms, $E$
is decreased.
Since there exists only a finite number of different
values $k_{ij}$ and the assignments are unique,
the algorithms converge in a finite number
of steps toward a fixed point $\vec w^*$ for which
$(\vec w^*)'=\vec w^*$ holds.

Consider the case of continuous $\vec w$.
Since $k_{ij}$ are discrete, $k_{ij}(\vec w)$
is constant in a vicinity of 
a fixed point $\vec w^*$ 
if no data points lie at the borders of 
two receptive fields. Then 
$E(\cdot)$ and $Q(\cdot,\vec w^*)$ are identical
in a neighborhood of $\vec w^*$ and thus,
a local optimum of $Q$ is also
a local optimum of $E$.
Therefore, if $\vec w$
can be varied in a real vector space, 
a local optimum of $E$ is found by the batch variant
if no data points are directly located at the
borders of receptive fields for the final solution. 

\section{Experiments}
We demonstrate the behavior of the algorithms
in different scenarios which cover a variety of characteristic situations.
All algorithms have been implemented based on the SOM Toolbox for Matlab 
\cite{somtoolbox}.
We used k-means, SOM,
Batch SOM, and NG with default
parameters as provided in the toolbox.
Batch NG and median versions of NG, SOM, and k-means have been implemented
according to the above formulas.
Note that, for all batch versions, prototypes which
lie at identical points of the data space do not separate in
consecutive runs.
Thus, the situation of exactly identical prototypes 
must be avoided.
For the euclidian versions, this situation is a set
of measure zero if prototypes are initialized at different positions.
For median versions, however, it can easily happen
that prototypes
become identical due to a limited number
of different positions in the data space,
in particular for small data sets.
Due to this fact, we add a small amount of noise to the distances
in each epoch in order to separate identical prototypes.
Vectorial training sets are normalized prior to training using
z-transformation.
Initialization of prototypes takes place using small random values.
The initial neighborhood rate for neural gas is $\lambda=n/2$, 
$n$ being the number of neurons, and it is 
multiplicatively decreased during training.
For Median
SOM, we restrict to square lattices of $n=\sqrt n \times\sqrt n $ neurons
and a rectangular neighborhood structure,
whereby $\sqrt{n}$ is rounded to the next integer.
Here the initial neighborhood rate is $\sqrt n /2$.

\begin{figure}[tb]
\begin{center}
\epsfxsize=10.6cm
\epsffile{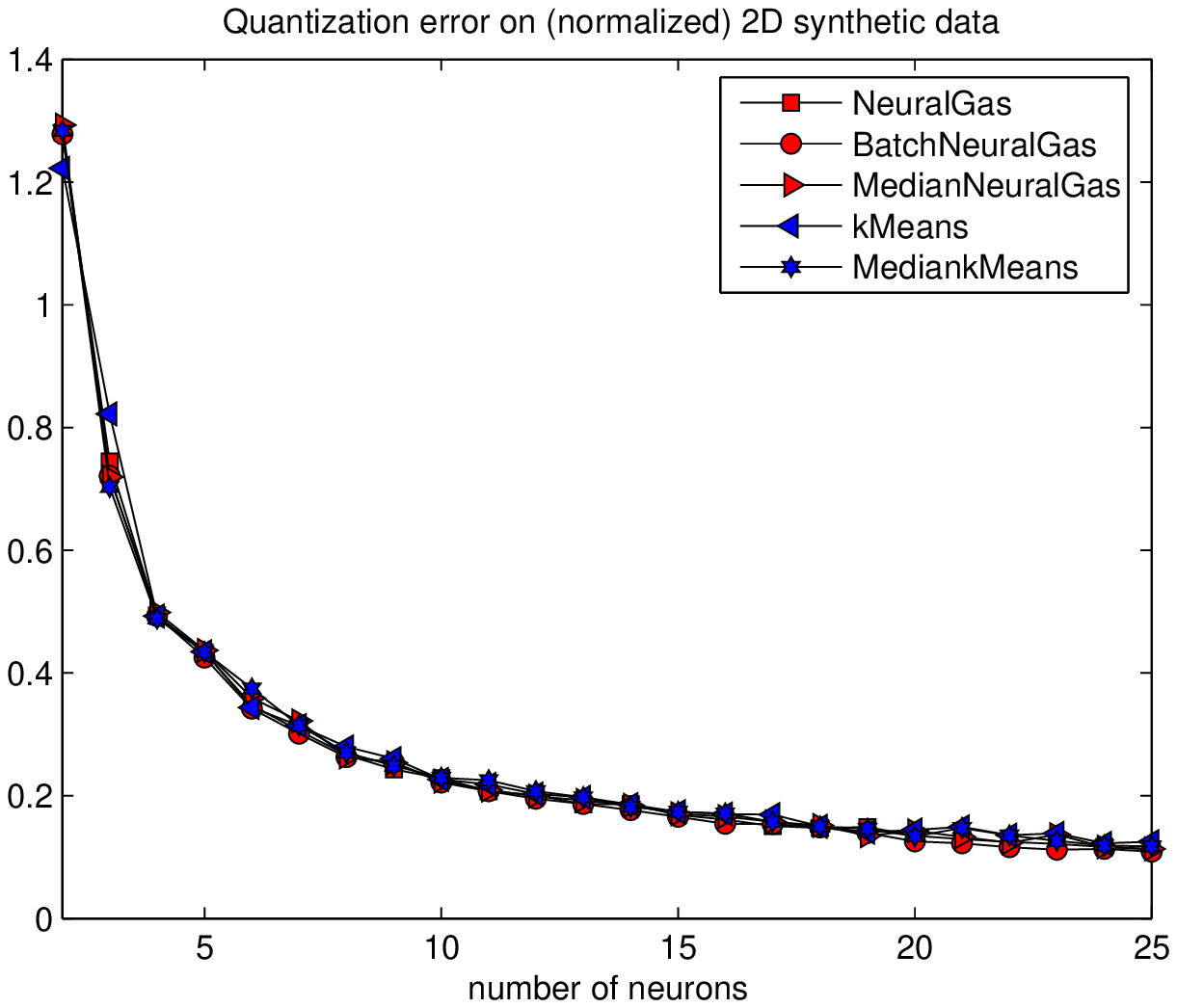}\\
\epsfxsize=10.6cm
\epsffile{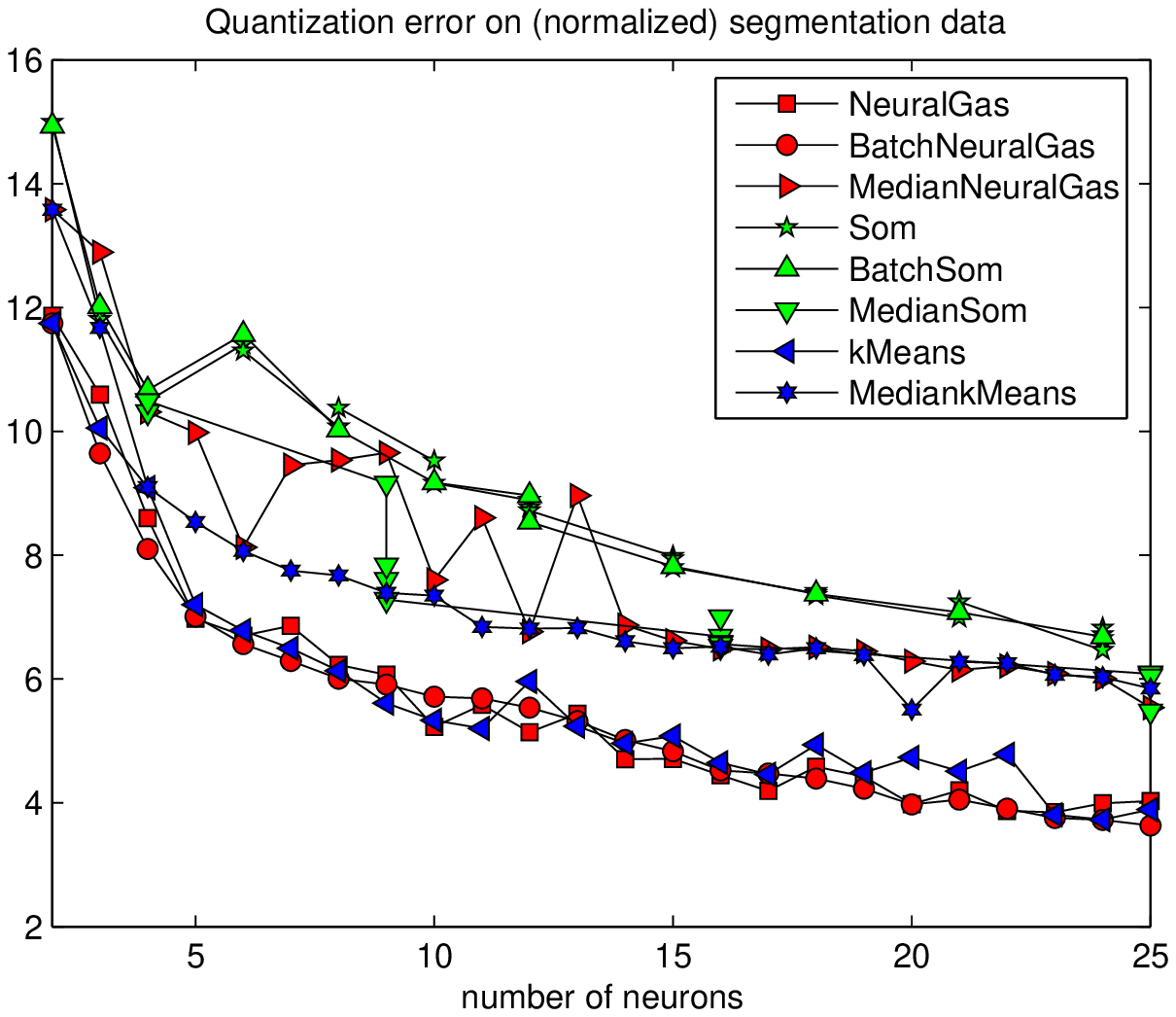}
\end{center}
\caption{Mean quantization error of the methods for the synthetic data set (top)
and the segmentation data set (bottom).}
\label{synth}
\end{figure}

\begin{figure}[tb]
\begin{center}
\epsfxsize=10.6cm
\epsffile{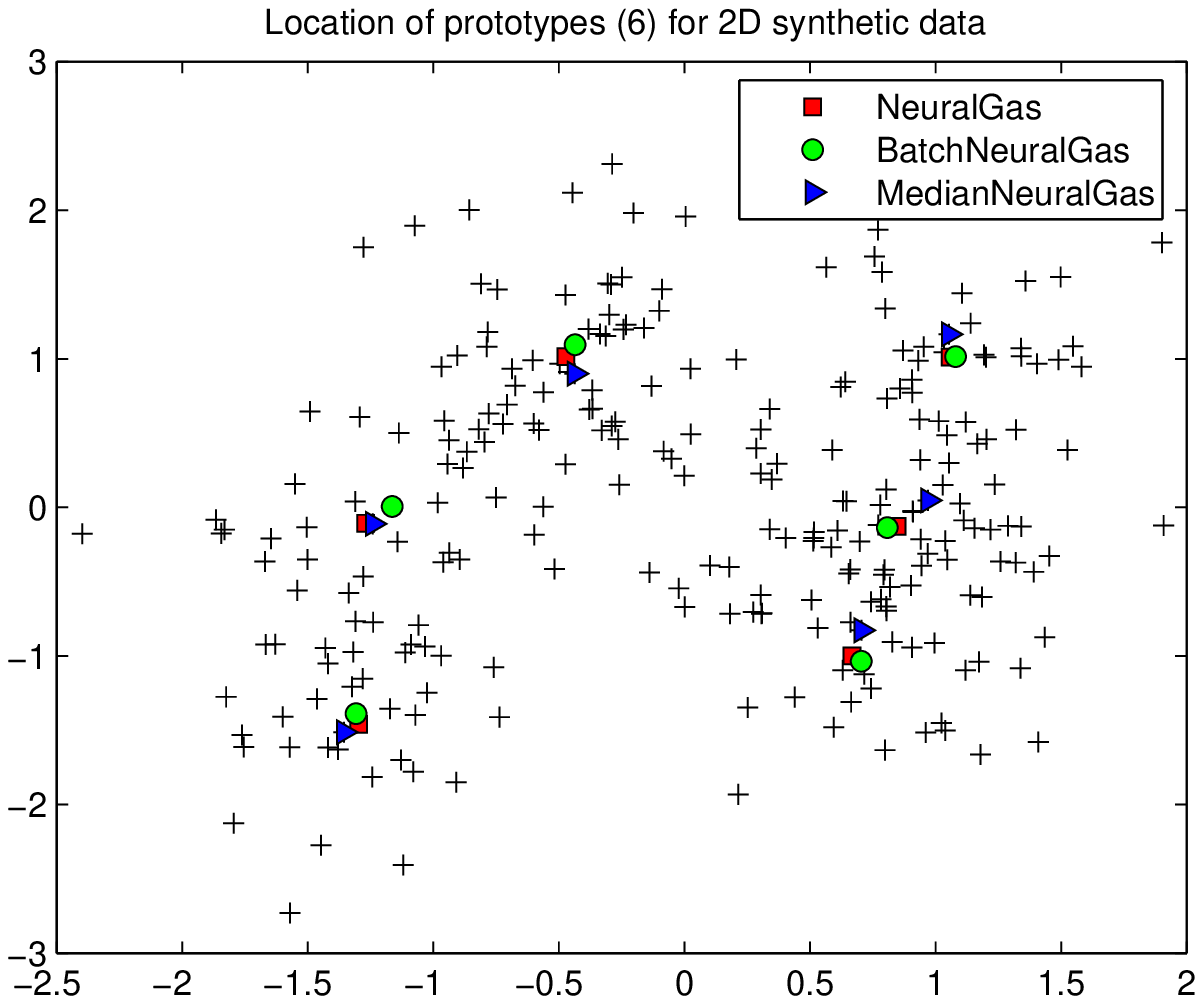}\\
\epsfxsize=10.6cm
\epsffile{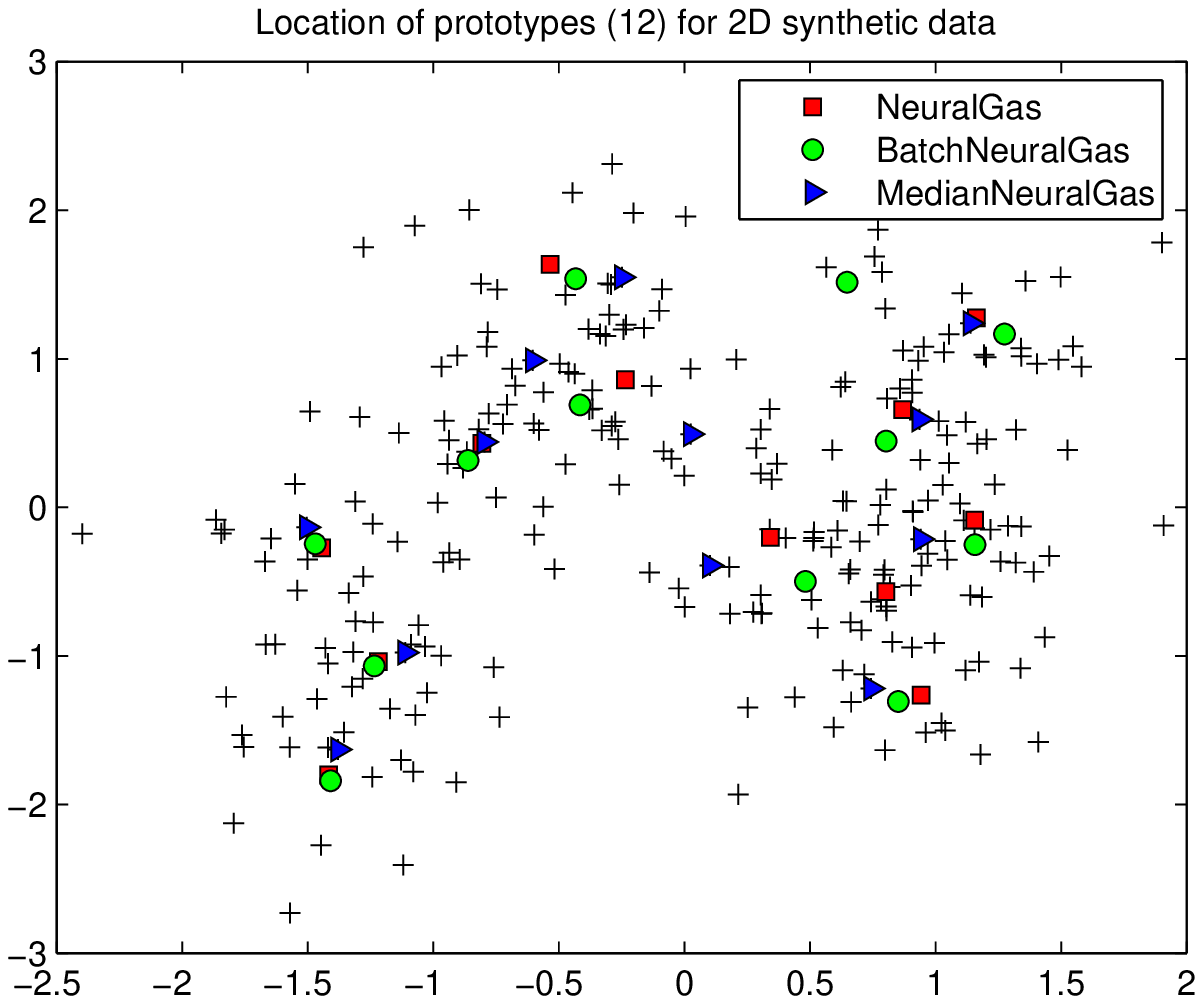}
\end{center}
\caption{Location of the prototypes for the synthetic data set for
different variants of NG.}
\label{seg}
\end{figure}

\subsection{Synthetic data}
The first data set is the two-dimensional synthetic data set
from \cite{ripley} consisting of
$250$ data points and $1000$ training points.
Clustering has been done using 
$n=2$, \ldots, $25$ prototypes, resp.\ the closest 
number of prototypes implemented by a rectangular lattice for SOM.
Training takes place for $5n$ epochs.

The mean quantization error $\sum_{i=1}^n\sum_{j=1}^p\chi_{I(\vec x^j)}(i)\cdot d(\vec x^j,\vec w^i)/p$ on the test set and 
the location of prototypes within the training set are depicted
in Figs.~\ref{synth} and \ref{seg}.
Obviously, the location of prototypes coincides for
different versions of NG.
This observation also holds for different numbers of prototypes, whereby
the result is subject to random fluctuations for larger numbers.
For k-means, idle prototypes can be observed for large $n$.
For Batch SOM and standard SOM, the quantization error is
worse (ranging from $1.7$ for  $2$ neurons up to $0.3$ for $24$ neurons,
not depicted in the diagram),
which can be attributed to the fact that the map does not fully unfold
upon the data set and edge effects remain,
which could be addressed to a small but nonvanishing
neighborhood in the convergent phase in standard implementations of SOM
which is necessary to preserve topological order.
Median SOM (which has been directly implemented in analogy to
Median NG)
yields a quantization error
competitive to NG.
Thus, Batch and Median NG allow to achieve results competitive to NG
in this case, however, using less effort.

\subsection{Segmentation data}
The segmentation data set from the UCI repository
consists of $210$ (training set) resp.\ $2100$ (test set) $19$ dimensional
data points which are obtained as pixels from outdoor images preprocessed
by standard filters such as averaging, saturation, intensity, etc.
The problem is interesting since
it contains high dimensional and only sparsely covered
data.
The quantization error obtained for the test set is
depicted in Fig.~\ref{synth}.
As beforehand, SOM suffers from the restriction of the topology.
Neural gas yields very robust
behavior, whereas for k-means, idle prototypes can be observed.
The median versions yield a larger quantization error compared
to the vector-based algorithms.
The reason lies in the fact that a high dimensional data set with only
few training patterns is considered, such that the
search space for median algorithms is small in these cases and
random effects and restrictions account for the increased error.

\begin{table}[b]
\begin{center}
\begin{tabular}{l|llllllll}
&NG& batch & median  & SOM & batch  & median  & kmeans & median\\[-0.7em]
&& NG &  NG &  &  SOM &  SOM &  &  kmeans\\\hline\hline
\multicolumn{9}{l}{quantization error}\\\hline
train&0.0043  & 0.0028  & 0.0043  & 0.0127  & 0.0126 &  0.0040& 0.0043 &  0.0046\\
test&0.0051  & {\bf 0.0033}  & 0.0048  & 0.0125  & 0.0124 &  0.0043&   0.0050 &  0.0052\\\hline\hline
\multicolumn{9}{l}{classification error}\\\hline
train&0.1032  & 0.0330  & 0.0338  & 0.2744  & 0.2770 &  0.0088&0.1136  & 0.0464\\
test&0.1207  & 0.0426  & 0.0473  & 0.2944  & 0.2926 &  {\bf 0.0111} & 0.1376 &  0.0606
\end{tabular}
\end{center}
\caption{Quantization error and
classification error for posterior labeling for
training and test set (both are of
size about $1800$). The mean over $5$ runs is reported.
The best results on the test set is depicted in boldface.}
\label{check}
\end{table}

\subsection{Checkerboard}
This data set is taken from \cite{npl}.
Two-dimensional data are arranged on a checkerboard, resulting in
$10$ times $10$ clusters, each consisting
of $15$ to $20$ points.
For each algorithm,
we train $5$ times $100$ epochs for $100$ prototypes.
Obviously, the problem is highly multimodal
and usually the algorithms  do not find all clusters.
The number of missed clusters can easily be 
judged in the following way:
the clusters are labeled consecutively using labels $1$ and $2$
according to the color black resp.\ white of the data
on the corresponding field of the checkerboard.
We can assign labels to
prototypes 
a posteriori based on a majority vote on the training set.
The number of errors which arise from this classification on an
independent
test set count the number of missed clusters, since
$1\%$ error roughly corresponds to one
missed cluster.

The results are collected in Tab.~\ref{check}.
The smallest quantization error is obtained by Batch NG, the smallest
classification error can be found for Median SOM.
As beforehand, the implementations for SOM and Batch SOM
do not fully unfold the map among the data.
In the same way online NG does not
achieve a small error because of a restricted number of epochs
and a large data set
which prevents online NG from full unfolding.
K-means also shows a quite high error (it misses more than $10$ clusters)
which can be explained by the existence of multiple
local optima in this setting, i.e.\ the sensitivity of k-means
with respect to initialization of prototypes.
In contrast, Batch NG and Median
NG find all but $3$ to $4$ clusters. Median SOM even finds
all but only $1$ or $2$ clusters since the topology
of the checkerboard exactly matches the
underlying data topology consisting of $10\times10$ clusters.
Surprisingly, also Median k-means shows quite good behavior,
unlike k-means itself, which might be due to
the fact that the generalized medians 
enforce the prototypes to settle within the clusters.
Thus, median versions and
neighborhood cooperation seem beneficial in this task due
to the multiple modes.
Batch versions show much better behavior than their online
correspondents, due to a faster convergence of the algorithms.
Here, SOM suffers from border effects, whereas Median SOM
settles within the data clusters, whereby the topology
mirrors precisely the data topology.
Both, Batch NG and Median NG, yield quite good
classification results which are even competitive
to supervised prototype-based
classification results as reported in \cite{npl}.

\begin{figure}[tb]
\begin{center}
\epsfxsize=9.0cm \epsffile{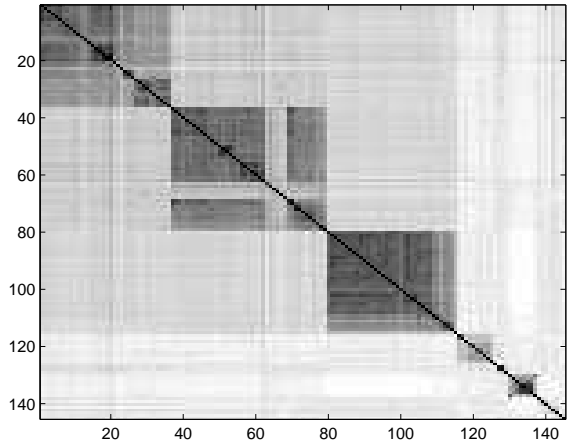}
\end{center}
\caption{Distance matrix for protein data.}
\label{matrix}
\end{figure}

\subsection{Proximity data -- protein clusters}
\begin{figure}[tb]
\begin{center}
\epsfxsize=9.0cm\epsffile{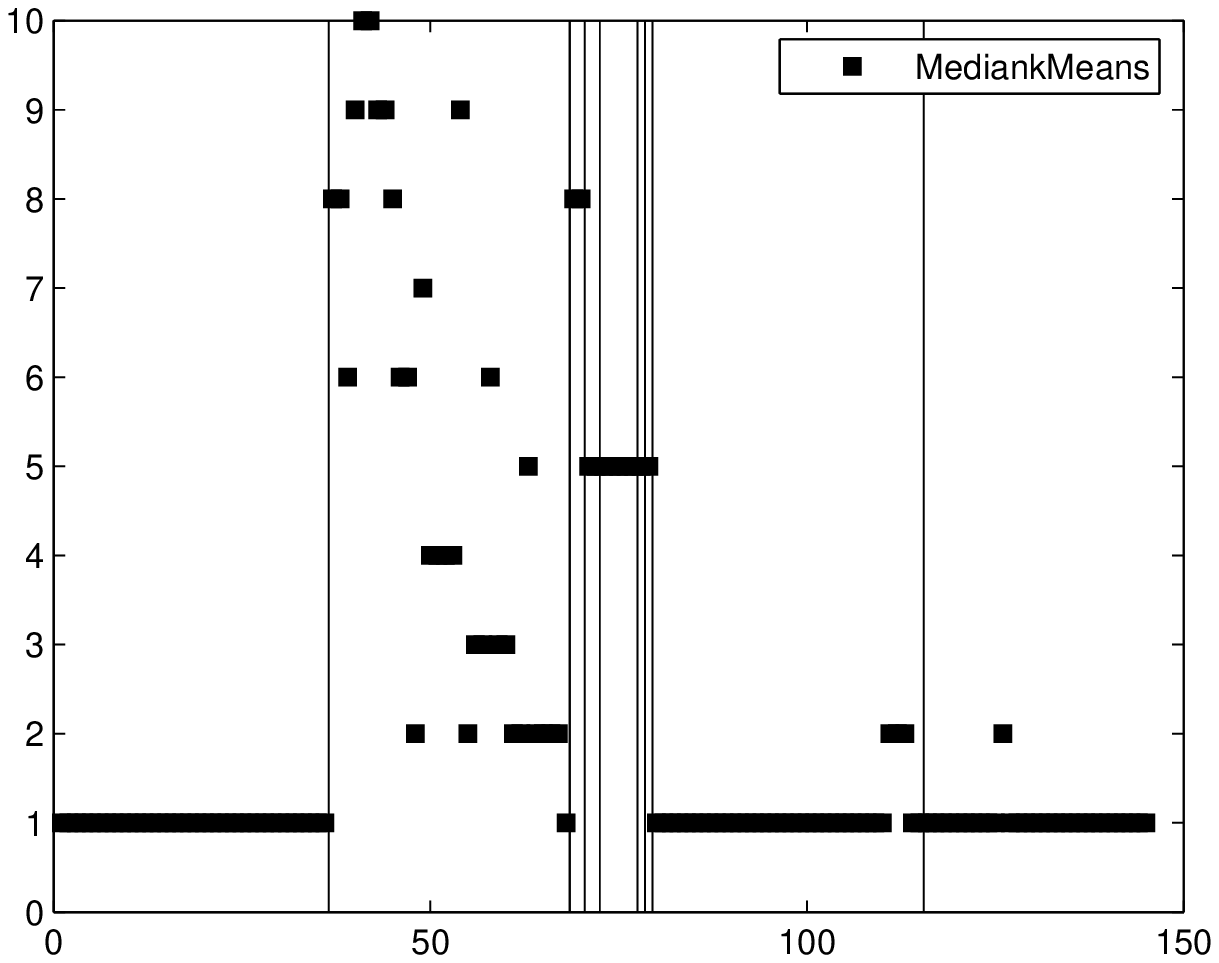}\\
\epsfxsize=9.0cm\epsffile{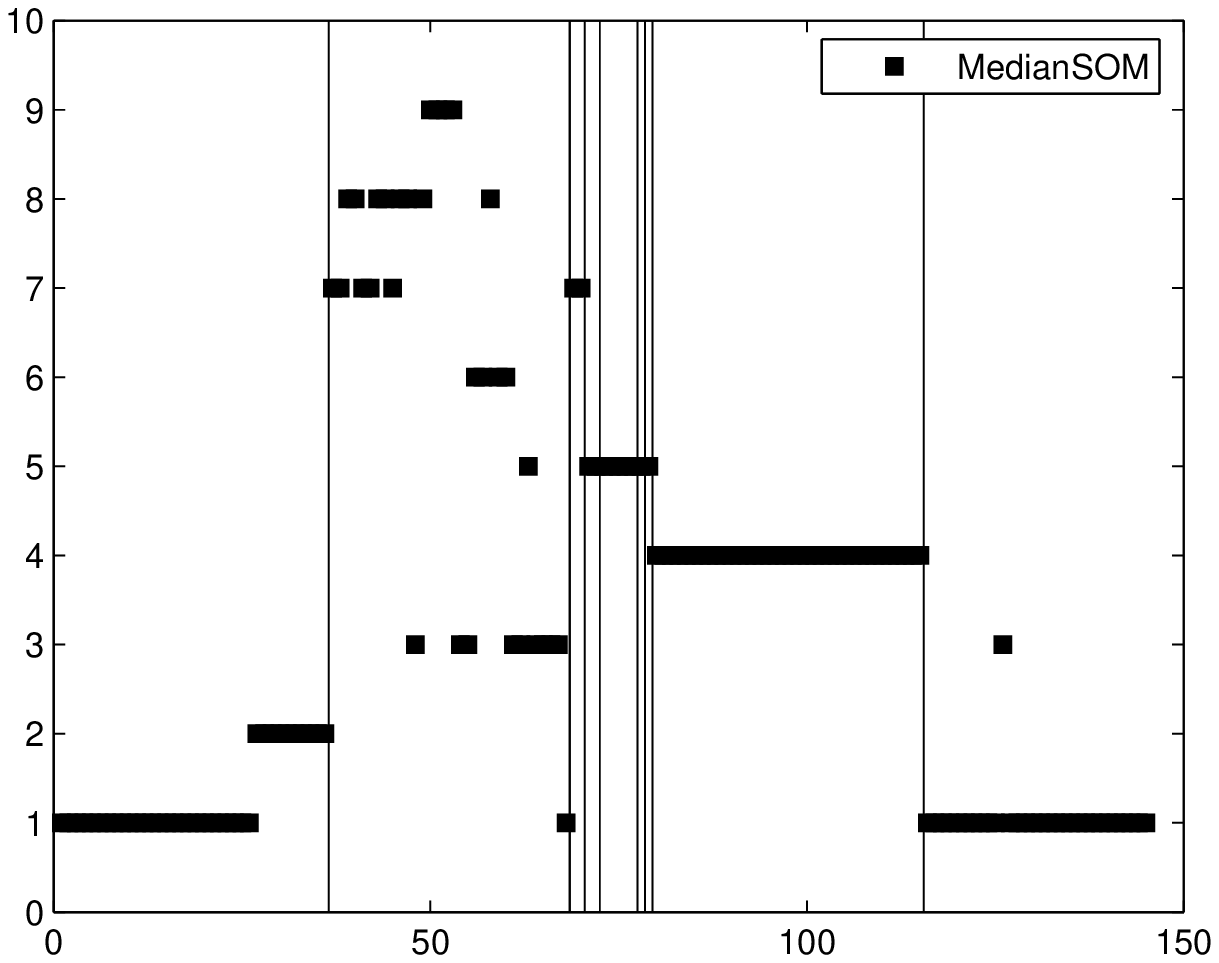}\\
\epsfxsize=9.0cm\epsffile{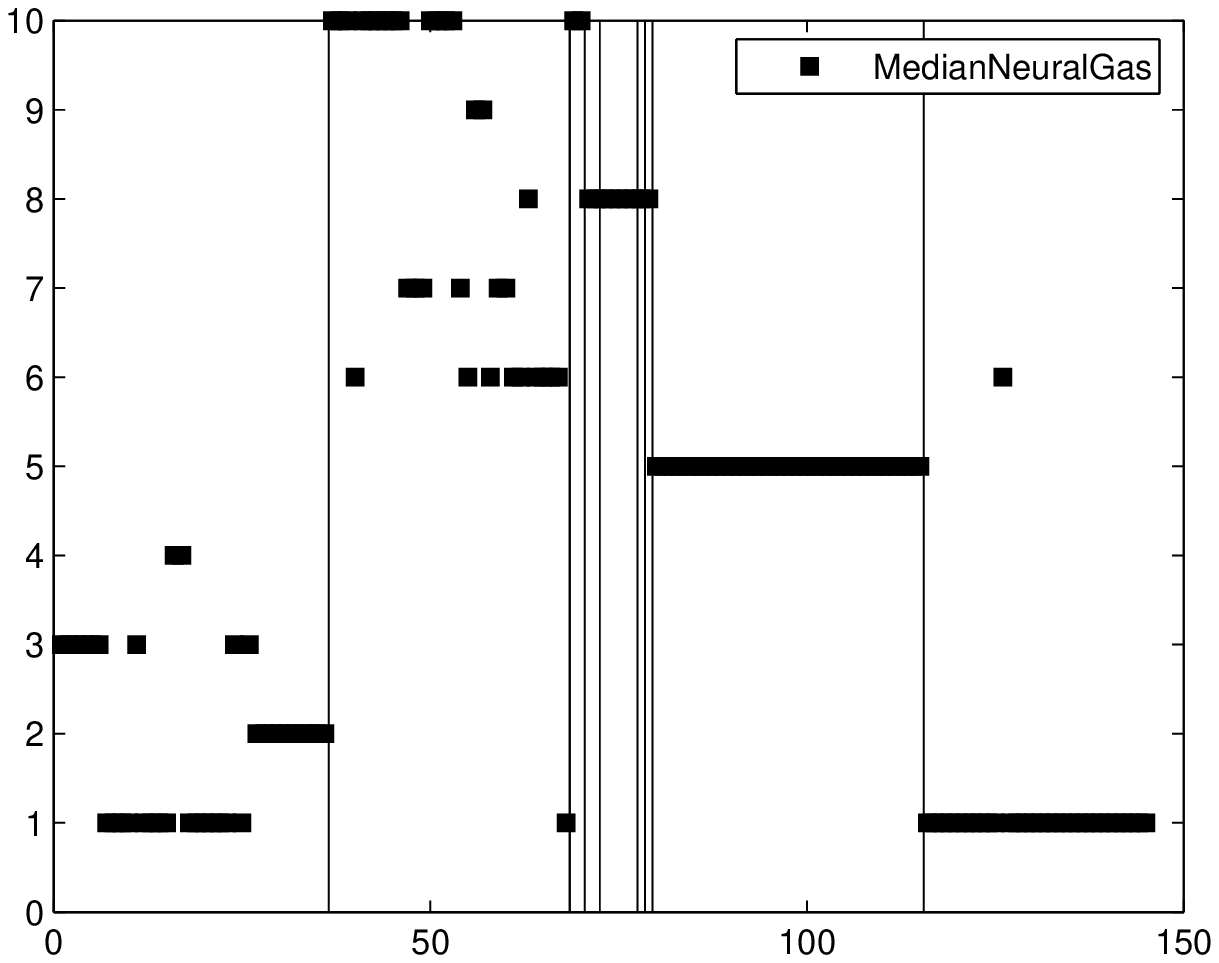}
\end{center}
\caption{
Typical results for median classification and $10$ prototypes.
The x-axes shows the protein number, the y-axes its winner neuron.
The vertical lines indicate an expert classification into different protein families (from left to right: hemoglobin $\alpha$, $\beta$, $\delta$,
$\epsilon$, $\gamma$, F, myoglobin, others).}
\label{matrixcluster}
\end{figure}

We used the protein data set described in 
\cite{protein} and \cite{seo}: the dissimilarity of $145$ globin
proteins of different families 
is given in matrix form as depicted in Fig.~\ref{matrix}.
Thereby, the matrix is
determined based on sequence alignment using biochemical and structural
information.
In addition, prior information about the underlying
protein families is available, i.e.\ a prior
clustering into semantically meaningful classes
of the proteins is known:
as depicted in Fig~\ref{matrixcluster} by
vertical lines, the first 42 proteins belong to hemoglobin $\alpha$,
the next clusters denote hemoglobin $\beta$, $\delta$, etc.
Thereby, several clusters are rather small, comprising
only few proteins (one or two).
In addition, the cluster depicted on the right has a very large
intercluster distance.

Since only a proximity matrix is available, we cannot apply standard
NG, k-means, or SOM, but we can rely on the median versions.
We train all three median versions $10$ times using $10$ prototypes and $500$ epochs.
The mean quantization errors (and variances) are
$3.7151$ ($0.0032$) for Median NG
$3.7236$ ($0.0026$) for Median SOM, and $4.5450$ ($0.0$)
for Median k-means, thus k-means yields worse results compared to NG and SOM
and neighborhood integration clearly seems beneficial in this 
application scenario.

We can check whether the decomposition into clusters by means of the prototypes
is meaningful by comparing the receptive fields of the ten prototypes
to the prior semantic clustering.
Typical results are depicted in Fig.~\ref{matrixcluster}.
The classification provided by experts is indicated by vertical lines in the 
images. The classification by the respective median method is
indicated by assigning a value on the y-achses to each pattern corresponding
to the number of its winner neuron (black squares in the figure).
Thus, an assignment of all or nearly all patterns in one
semantic cluster to one or few dedicated prototypes gives a hint for the fact
that median clustering finds semantically meaningful entities.

All methods detect the first cluster (hemoglobin $\alpha$) and
neural gas and SOM also detect the eighth cluster (myoglobin).
In addition, SOM and NG group together elements of clusters 
two to seven in a reasonable way.
Thereby, according to the variance in the clusters, more than one
prototype is used for large clusters and small clusters containing only
one or two patterns are grouped together. 
The elements of the last two clusters
have a large intercluster distance
such that they are grouped together into some (random) cluster
for all methods.
Note that the goal of NG and SOM
is a minimization of
their underlying cost function, such that the 
cluster border can lie between semantic
clusters for these methods.
Thus, the results obtained by SOM and NG are reasonable and they detect 
several semantically meaningful clusters.
The formation of relevant clusters is also supported when training
with a different number of prototypes

\subsection{Proximity data -- chicken pieces silhouettes}
The data set as given in \cite{acv} consists of silhouettes of 
$446$ chicken pieces of different classes including wings, backs, drumsticks, thighs,
and breasts.
The task is a classification of the images (whereby the silhouettes are not
oriented) into the correct class.
As described in \cite{spillmann}, a preprocessing of the images resulting in
a proximity matrix can cope with the relevant properties of
the silhouette and rotation symmetry:
the surrounding edges are detected and discretized into small
consecutive line segments of $20$ pixels per segment. 
The images are then represented by the 
differences of the angles of consecutive line segments.
Distance computation takes place as described in \cite{bunke}
by a rotation and mirror symmetric variant of the edit distance
of two sequences of angles, whereby
the costs for a substitution of two angles is given
by their absolute distance, the costs for deletion and insertion
are given by $k=60$.

We train Median k-means, Median NG, and Median SOM with different
numbers of neurons for $500$ epochs, thereby annealing the neighborhood
as beforehand.
The results on a training and test set of the same size, averaged over ten runs,
are depicted in Tab.~\ref{chicken_pieces}.
Obviously, a posterior labeling of prototypes obtained
by median clustering allows to achieve a classification
accuracy of more than $80\%$.
Thereby, overfitting can be
observed for all methods due to the large number of
prototypes compared to the training set ($50$ neurons constitute
about $1/4$th of the training set!).
However, Median NG and Median SOM are less prone to
this effect due to their inherent regularization given by the
neighborhood integration.

\begin{table}[bt]
\begin{center}
\begin{tabular}{l|ll|ll|ll}
neurons&\multicolumn{2}{c|}{Median k-means}&\multicolumn{2}{c|}{Median NG}&\multicolumn{2}{c}{Median SOM}\\
&train&test&train&test&train&test\\\hline
10&0.54&0.52&0.57&{\bf 0.61}&0.53&0.59\\
20&0.71&0.61&0.67&{\bf 0.65}&0.61&0.61\\
30&0.77&0.63&0.73&{\bf 0.64}&0.69&0.61\\
40&0.85&{\bf 0.79}&0.80&0.75&0.74&0.69\\
50&0.90&0.79&0.84&{\bf 0.80}&0.76&0.68\\
60&0.88&0.82&0.88&{\bf 0.83}&0.80&0.73\\
70&0.93&0.82&0.89&{\bf 0.84}&0.89&0.78\\
80&0.94&0.82&0.92&{\bf 0.84}&0.87&0.78\\
90&0.95&0.81&0.93&{\bf 0.84}&0.87&0.78\\
100&0.96&{\bf 0.83}&0.94&{\bf 0.83}&0.88&0.80
\end{tabular}
\end{center}
\caption{Results for the median variants for different numbers
of neurons on the chicken-piece-silhouettes data base.
The best test classifications are depicted in bold.}\label{chicken_pieces}
\end{table}

\subsection{Proximity data -- chromosomes}
The Copenhagen Chromosomes Database \cite{lpg}
consists of $4400$ descriptions of chromosomes by their silhouettes
in images. A chromosome is described by a sequence
over the alphabet $\{1,\ldots,6\}$,
whereby the number describes the thickness of the density profile
of the protein at the corresponding position.
The difference between two profiles is determined by alignment
assigning the costs $|x-y|$ to substitutions of $x$ and $y$,
and assigning the costs $4.5$ to insertions and deletions,
as described in \cite{spillmann}.
There are $21$ different classes. The set is divided into a training
and test set of the same size.

We train median clustering with
different numbers of neurons and $100$ cycles. The classification accuracy
on a training and test set, averaged over $10$ runs, is depicted in
Tab.~\ref{chromosome}.
As beforehand, a classification accuracy of $80\%$ can be achieved.
Thereby, Median NG shows the best results on the
test set for almost all numbers of neurons,
accompanied by a good generalization error due
to the inherent regularization by means of neighborhood cooperation.

\begin{table}[bt]
\begin{center}
\begin{tabular}{l|ll|ll|ll}
neurons&\multicolumn{2}{c|}{Median k-means}&\multicolumn{2}{c|}{Median NG}&\multicolumn{2}{c}{Median SOM}\\
&train&test&train&test&train&test\\\hline
10&0.31&0.25&0.43&{\bf 0.40}&0.40&0.34\\
20&0.52&{\bf 0.45}&0.46&0.42&0.54&0.52\\
30&0.64&0.57&0.70&{\bf 0.66}&0.57&0.53\\
40&0.75&0.62&0.75&{\bf 0.71}&0.69&0.63\\
50&0.78&0.73&0.79&{\bf 0.74}&0.75&0.67\\
60&0.80&0.74&0.83&{\bf 0.78}&0.75&0.67\\
70&0.75&0.68&0.82&{\bf 0.77}&0.69&0.60\\
80&0.82&0.75&0.83&{\bf 0.78}&0.68&0.58\\
90&0.82&0.74&0.82&{\bf 0.76}&0.73&0.65\\
100&0.82&0.76&0.86&{\bf 0.81}&0.78&0.72
\end{tabular}
\end{center}
\caption{Classification accuracy on the chromosome data set
for different numbers of neurons.
The best results on the test set are depicted in bold.}
\label{chromosome}
\end{table}

\section{Conclusions}
We have proposed Batch NG derived from the NG cost function
which allows fast training for a priorly given data set.
We have shown that the method converges and it optimizes the same cost function
as NG by means of a Newton method. 
In addition, the batch formulation opens the way towards
general proximity data by means of the generalized median.
These theoretical discussions were supported by experiments for
different
vectorial data where the results of Batch NG and NG are very similar.
In all settings, the quality of Batch NG was at least competitive to
standard NG, whereby training takes place in a fraction of the time
especially for high-dimensional input data due to the radically reduced
number of updates of a prototype.
Unlike k-means, NG is not sensitive to initialization and,
unlike SOM, it automatically determines a data optimum lattice,
such that a small quantization error can be achieved and topological initialization is not crucial.

Median NG restricts the adaptation to locations within the data set
such that it can be applied to non-vectorial data.
We compared Median NG to its alternatives for vectorial data
observing that competitive results arise if enough data are available.
We added several experiments including proximity data where we could obtain
semantically meaningful grouping as demonstrated by a
comparison to known clusters resp.\ a validation of the
classification error when used in conjunction with posterior labeling.
Unlike SOM,
NG solely aims at data clustering and not data visualization, such that
it can use a data optimum lattice and it is not restricted
by topological constraints. Therefore better
results can often be obtained in terms of the quantization error
or classification. If a visualization of the output of
NG is desired, a
subsequent visualization of the prototype vectors is possible
using fast standard methods for
the reduced set of prototypes such as multidimensional scaling \cite{mds}.

Thus, very promising results could be achieved which have been accompanied
by mathematical guarantees for the convergence of the algorithms.
Nevertheless, several issues remain:
for sparsely covered data sets, median versions might not have
enough flexibility to position the prototypes since only
few locations in the data space are available.
We have already demonstrated this
effect by a comparison of batch clustering
to standard euclidian clustering in such a situation.
It might be worth investigating metric-specific possibilities
to extend the adaptation space for the prototypes in such
situations, as possible e.g.\ for the edit distance,
as demonstrated in \cite{guenter} and \cite{somervuo}.

A problem of Median NG is given by the complexity of one 
cycle, which is quadratic in the number of patterns.
Since optimization of the exact computation as proposed in \cite{fabrice}
is not possible, heuristic variants which restrict the computation
to regions close to the winner seem particularly promising because
they have a minor effect on the outcome.
A thorough investigation of the effects of such restriction
will be investigated both theoretically and experimentally in future
work.

Often, an appropriate metric or proximity matrix is not fully
known a priori. The technique of learning metrics, which
has been developed for both, supervised as well as
unsupervised prototype-based methods \cite{npl,kaski}
allows a principled integration of secondary knowledge into the framework
and adapts the metric accordingly, thus getting around the
often problematic \lq garbage-in-garbage-out\rq\ problem of
metric-based approaches. It would be interesting to
investigate the possibility to enhance median versions for
proximity data by an automatic adaptation of the distance matrix during training
driven by secondary information.
A recent possibility to combine vector quantizers with prior
(potentially fuzzy) label information has been proposed in \cite{flng}
by means of a straightforward extension of the underlying
cost function of NG. This approach can immediately be transferred to
a median computation scheme since a well-defined cost function
is available, thus opening the way towards supervised 
prototype-based median fuzzy classification for non-vectorial data.
A visualization driven by secondary label information can be developed
within the same framework substituting the irregular NG lattice by a SOM
neighborhood and incorporating Heskes' cost function.
An experimental evaluation of this framework is the subject
of ongoing work.

\end{document}